\documentclass[a4paper,12pt,leqno]{article}
\usepackage{latexsym}
\usepackage[all]{xy}

\usepackage{amssymb} 
\usepackage{amsmath} 
\usepackage{theorem}

\usepackage{theorem}
\usepackage{amscd}
\usepackage{latexsym}
\usepackage{enumerate}
\usepackage{ascmac,fancybox}
\usepackage[dvips]{graphicx,color}

\def\bfk{{\mathbf{k}}}
\def\bfm{{\mathbf{m}}}
\def\bfg{{\mathbf{g}}}
\def\Z{{\mathbb{Z}}}
\def\R{{\mathbb{R}}}
\def\A{{\mathcal{A}}}
\def\G{{\mathcal{G}}}

\DeclareMathOperator{\Der}{Der}

\DeclareMathOperator{\Poin}{Poin}

\numberwithin{equation}{section}

\newcommand{\owari}{\hfill$\square$}

\newtheorem{theorem}{Theorem}[section]
\newtheorem{prop}[theorem]{Proposition}
\newtheorem{cor}[theorem]{Corollary}
\newtheorem{lemma}[theorem]{Lemma}
\newtheorem{define}[theorem]{Definition}

\newtheorem{example}[theorem]{Example}

\title{
Equivariant multiplicities 
of Coxeter arrangements and 
invariant bases}
\author{Takuro Abe
\footnote
{
Supported by JSPS Grants-in-Aid for Young Scientists
(B)
No. 21740014.
Department of Mechanical Engineering and Science,
Kyoto University,
Kyoto 606-8501, Japan.
email:abe.takuro.4c@kyoto-u.ac.jp
}
\,
Hiroaki Terao
\footnote
{
Supported by JSPS Grants-in-Aid, Scientific Research
(B) 
No. 21340001.
Department of Mathematics, Hokkaido University, 
Sapporo, Hokkaido 060-0810, Japan.
email:terao@math.sci.hokudai.ac.jp
}
\,
Atsushi Wakamiko
\footnote
{{2-8-11
Aihara, Midori-ku,
Sagamihara-shi, Kanagawa,
252-0141
Japan.}
email:atsushi.wakamiko@gmail.com
}}

\begin{document}

\maketitle

\begin{abstract}
Let $\A$ be an irreducible Coxeter arrangement and $W$ be
its Coxeter group. Then $W$ naturally acts on $\A$. 
A multiplicity $\bfm : \A\rightarrow \Z$ is said to be 
equivariant when $\bfm$ is constant on each $W$-orbit of 
$\A$. 
In this article, we prove that the 
multi-derivation
module
$D(\A, \bfm)$ is a free module
whenever $\bfm$ is equivariant
by explicitly constructing a basis,
which generalizes the main theorem of \cite{T02}. 
The main tool is a primitive derivation
and its covariant derivative.
Moreover, 
we show that
the $W$-invariant part
$D(\A, \bfm)^{W}$
for any multiplicity $\bfm$
is a free module over 
the $W$-invariant subring.
\end{abstract}

\section{Introduction}
\label{introduction} 
Let $V$ be an $\ell$-dimensional Euclidean space
 with 
an inner product $I: V \times V \rightarrow \R$.
Let
$S$
denote the symmetric algebra of the dual space
$V^{*} $ and $F$ 
be its quotient field.
Let
$\Der_{S}$ be the $S$-module of $\R$-linear derivations 
from $S$ to itself.
Let 
$\Omega^1_S$ be the $S$-module of regular $1$-forms. 
Similarly define $\Der_{F} $ and $\Omega^{1}_{F}$
over $F$.
The dual inner product $I^{*} : V^{*} \times V^{*} 
\rightarrow \R$ naturally induces an $F$-bilinear
form   $I^{*} : \Omega^{1}_{F}  \times \Omega^{1}_{F}
\rightarrow F$. 
Then one has an $F$-linear bijection
\[
I^{*} : \Omega_{F}^{1} \rightarrow \Der_{F}  
\]
defined by
$\left[I^{*} (\omega)\right](f)
:=
I^{*}(\omega, df) $ 
for $f\in F$.  

Let $\A$ be an irreducible Coxeter arrangement
with its Coxeter group $W$. 
For each $H\in\A$, choose $\alpha_{H} \in V^{*} $
with $H = \ker(\alpha_{H})$.
Let $Q = \prod_{H\in\A} \alpha_{H} \in S$.
Recall 
the $S$-module of logarithmic forms 
\begin{eqnarray*}
\Omega^1(\A,\infty) =
\{
\omega \in \Omega^1_F &|&
Q^{N}\omega 
\mbox{~and~}
(Q/\alpha_{H})^{N} \omega \wedge d \alpha_H
\mbox{~are both regular}\\
&&\mbox{for any~}H\in \A
\mbox{~and~} N \gg 0\}
\end{eqnarray*}
and the $S$-module of 
logarithmic derivations
\begin{eqnarray*}
D(\A, -\infty)&=&I^{*} (\Omega^{1}(\A, \infty))
\end{eqnarray*}
from \cite{AT08}. 
A map $\bfm: \A \rightarrow \Z$ is called a multiplicity.
For an arbitrary multiplicity, let 
\begin{eqnarray*}
D(\A,\bfm)&=&\{\theta \in D(\A, -\infty)
\mid \theta(\alpha_H) \in \alpha_H^{\bfm(H)}\ S_{(\alpha_{H})}  
\mbox{ for all }H \in \A\},\\
\Omega^1(\A,\bfm) &=&
(I^{*})^{-1} D(\A,-\bfm),
\end{eqnarray*}
where 
$S_{(\alpha_H)}$ is the localization of $S$ at the prime ideal
$(\alpha_{H})$. 
{
These two modules were introduced in \cite{S80} (when
 $\bfm$ is constantly equal to one),
in \cite{Z89} (when
${\rm im}(\bfm)\subset \Z_{>0}$),
and in \cite{A08, AT08, AT09} (when
$\bfm$ is arbitrary).}
A derivation $0\neq\theta\in\Der_{F} $ is said to be {\bf
homogeneous
of degree $d$}, or $\deg\theta=d$,  if $\theta(\alpha)\in F$ is homogeneous of degree $d$ 
whenever $\theta(\alpha)\neq 0$ ($\alpha\in V^{*}$).
A multiarrangement 
$(\A,\bfm)$ is called to be {\bf free}
 with {\bf exponents}
 $\exp(\A,\bfm)=(d_1,\ldots,d_\ell)$ 
if $D(\A,\bfm)=\oplus_{i=1}^\ell S \cdot \theta_i$ 
with a homogeneous basis $\theta_i$ such that $\deg (\theta_i)=d_i\ (i=1,
\ldots,\ell)$. 
A multiplicity $\bfm : \A\rightarrow \Z$ is said to be 
{\bf equivariant} when $\bfm(H) = \bfm(wH)$ for any
$H\in\A$ and any $w\in W$, i.e.,
$\bfm$ is constant on each orbit. 
In this article we prove

\begin{theorem}
For any irreducible Coxeter arrangement $\A$ and any
equivariant multiplicity $\bfm$,
 the multiarrangement
$(\A, \bfm)$ is free.
\label{basis}
\end{theorem}

%
{
For a fixed arrangement
$\A$, we say that a multiplicity
$\bfm$ is {\bf 
free} if $(\A, \bfm)$ is free.  
Although we have a limited knowledge about
the set of all free multiplicities for a fixed 
irreducible Coxeter arrangement
$\A$,
it is known that there exist infinitely many non-free multiplicities
unless $\A$
 is either
one- or
two-dimensional
\cite{ATY}.  
Theorem \ref{basis} claims that any equivariant
multiplicity is free for any irreducible Coxeter arrangement.
}

When the $W$-action on $\A$ is transitive,
an equivariant multiplicity
 is
constant and
a basis was constructed in \cite{ST98,
T02, AY07, AT08}.
So we may assume,
in order to prove Theorem \ref{basis},  that
the $W$-action on $\A$ is not transitive.
In other words, we may only study the cases when $\A$ 
is of the
type
either
 $B_\ell,
F_4, G_2$ or 
$I_{2}(2n)\ (n\ge 4)$.
In these cases,
$\A$ has
exactly two $W$-orbits:
$\A = \A_{1} \cup \A_{2} $.
The orbit decompositions are explicitly
given by:
$B_{\ell}=
A_1^\ell
\cup 
D_{\ell} $,
$F_{4}
=
D_4 \cup D_4$, 
$G_{2} 
=
A_2 \cup A_2$ or 
$
I_{2}(2n) 
=
I_{2}(n) \cup I_{2}(n)
\ (n\ge 4)$.
Note that 
$A_1^\ell
$ 
is not irreducible.

When $\A$ is irreducible, the {\bf
primitive derivations} 
play the central role to define the
Hodge filtration introduced
by K. Saito.  (See \cite{S03} for example.)
For
$R:=S^{W}$,
let $D$ 
be an element of 
the lowest degree in $\Der_{R}$,
which is called
a primitive derivation corresponding to 
$\A$.
Then $D$ is unique up to a nonzero
constant multiple.
A theory of primitive
derivations
in the case of 
non-irreducible Coxeter arrangements
was introduced in 
\cite{AT09}.
Thus
we may consider a primitive derivation 
$D_{i}$ corresponding
with the orbit $\A_{i}$ 
$(1\le i \le 2)$. 
We only use $D_{1}$ because of symmetricity.
{
Note that
$D_{1} $ is
not unique up to a nonzero multiple
when $\A_{1} = A_{1}^{\ell}$
(non-irreducible). 
}
Denote the reflection groups of $\A_{i} $ by 
$W_{i}  $
 $(i=1,2)$.
The Coxeter arrangements 
 $B_\ell,
F_4, G_2$ and
$I_{2}(2n)\ (n\ge 4)$
are classified into two 
cases, that is,
(1)
the primitive derivation
$D_{1} $ 
{
can be
chosen to be}
$W$-invariant for
 $B_\ell$
and
$F_4$
(the first case)
while
(2)
$D_{1} $ 
is
$W_{2} $-antiinvariant for
$G_2$ and
$I_{2}(2n)\ (n\ge 4)$
(the second case)
as we will see in Section
4.
Since the second cases
are two-dimensional,
Theorem \ref{basis} holds true.
Thus the first case is the only
remaining case to prove Theorem \ref{basis}.
%

Let
\begin{align*} 
\nabla : 
&\Der_{F} \times \Der_{F} 
\longrightarrow
\Der_{F} \\
&~~~~~(\theta, \delta)
~~~~~~~\mapsto
~~\nabla_{\theta}\,\delta 
\end{align*} 
be the {\bf
Levi-Civita connection}
 with respect to 
the inner product $I$ on $V$. 
We need the following theorem for
 our proof of
Theorem \ref{basis}:

\begin{theorem}
\label{automorphism}
(\cite{AT08,AT09})
Let 
$D(\A, -\infty)^{W} $
 be the $W$-invariant part of 
 $D(\A, -\infty)$. Then
\[
\nabla_{D} : 
D(\A, -\infty)^{W} 
\overset{\sim}{\longrightarrow}
D(\A, -\infty)^{W} 
\]
is a $T$-linear automorphism
where
$
T := \{f \in R \mid Df=0\}.$
When $\A
=
\A_{1} \cup \A_{2} $ is the
orbit decomposition,
\[
\nabla_{D_{1} } : 
D(\A_{1}, -\infty)^{W_{1} } 
\overset{\sim}{\longrightarrow}
D(\A_{1}, -\infty)^{W_{1} } 
\]
is a $T_{1} $-linear automorphism
where
\[
R_{1} :=S^{W_{1} }, \,\,
T_{1}  := \{f \in R_{1} \mid D_{1} f=0\}.
\]
\end{theorem}

Let $E$ be the {\bf
Euler derivation
}
characterized by the equality
$E(\alpha) = \alpha$ 
for every $\alpha\in V^{*} $.  
Suppose that $\A = \A_{1} \cup \A_{2} $ is the orbit decomposition
and that the primitive derivation $D_{1} $ is $W$-invariant. 
Define
\[
E^{(p, q)}
:=
\nabla_{D}^{-q} \nabla_{D_{1} }^{q-p} E   
\]
 for $p, q\in\Z$.
Here,
thanks to Theorem \ref{automorphism},
 we may
interpret 
$\nabla_{D}^{m} = (\nabla_{D}^{-1})^{-m}  $ 
and
$\nabla_{D_{1} }^{m} = (\nabla_{D_{1} }^{-1})^{-m}  $ 
when  $m$ is negative.
Denote the equivariant multiplicity
$\bfm$ by $(m_{1}, m_{2})$ when 
$\bfm(H) = m_{1} \
(H\in \A_{1})$ 
and
$\bfm(H) = m_{2} \
(H\in \A_{2})$. 
Let $x_{1}, \dots , x_{\ell}$ be a basis for $V^{*} $ 
and
$P_1,\ldots,P_\ell$ be 
a set of basic invariants with respect to 
$W$: 
$R = \R[P_{1}, \dots , P_{\ell}].$
Let
$P^{(i)}_1,\ldots,P^{(i)}_\ell$ be 
a set of basic invariants with respect to 
$W_{i}$:
$R_{i}  = \R[P^{(i)} _{1}, \dots , P^{(i)}_{\ell}]$
 $(i=1,2)$.
Define
\[
d_{j} := \deg P_{j},
\,\,
d^{(i)}_{j} := \deg P^{(i)}_{j}
\,\,
(i = 1,2, 1\leq j\leq \ell).
\]
We assume
\[
d_{1} \leq 
d_{2}  \leq 
\dots
\leq
d_{\ell},
\,\,\,
d^{(i)}_{1} \leq 
d^{(i)}_{2}  \leq 
\dots
\leq
d^{(i)}_{\ell} 
\,\,\,(i = 1,2).
\]
Then $h := d_{\ell} $ is called the 
Coxeter number of $W$.  We call 
$h_{i} := \deg P^{(i)}_{\ell}$ 
the Coxeter number of $W_{i} $  
$(i=1, 2)$. 
We use the notation
\[
\partial_{x_{j}} :=
\partial /\partial x_{j},\, \,
\partial_{P_{j}} :=
\partial /\partial P_{j},\, \,
\partial_{P^{(i)}_{j}} :=
\partial /\partial P^{(i)}_{j}\, \, \, \,
(1\le j\le \ell, 1\le i\le 2).
\]
The following theorem gives an explicit construction of 
a basis:

\begin{theorem}
\label{basisconstruction} 
Let $\A$ be an irreducible Coxeter arrangement.
Suppose that $\A = \A_{1} \cup \A_{2} $ is the orbit decomposition
and that the primitive derivation $D_{1} $ 
is $W$-invariant. Then

\begin{itemize}
\item[(1)]
the $S$-module $D(\A,(2p-1,2q-1))$ is free with 
$W$-invariant basis 
$$
\nabla_{\partial_{P_1}} E^{(p,q)},\ldots,
\nabla_{\partial_{P_\ell}} E^{(p,q)}
$$
with $\deg \nabla_{\partial_{P_i}} E^{(p,q)} = p h_{1} + q h_{2} - d_{i} + 1$ 
for $i = 1, \dots, \ell$, 
\item[(2)]
the $S$-module $D(\A,(2p-1,2q))$ is free with basis 
$$
\nabla_{\partial_{P^{(1)}_1}} E^{(p,q)},\ldots,
\nabla_{\partial_{P^{(1)}_\ell}} E^{(p,q)}
$$
with $\deg \nabla_{\partial_{P^{(1)} _i}} E^{(p,q)} = p h_{1} + q h_{2} - 
d^{(1)}_{i} + 1$ 
for $i = 1, \dots, \ell$, 

\item[(3)]
the $S$-module $D(\A,(2p,2q-1))$ is free with basis 
$$
\nabla_{\partial_{P^{(2)}_1}} E^{(p,q)},\ldots,
\nabla_{\partial_{P^{(2)}_\ell}} E^{(p,q)}
$$
with $\deg \nabla_{\partial_{P^{(2)} _i}} E^{(p,q)} = p h_{1} + q h_{2} - 
d^{(2)}_{i} + 1$ 
for $i = 1, \dots, \ell$, 

\item[(4)]
the $S$-module $D(\A,(2p,2q))$ is free with basis 
$$
\nabla_{\partial_{x_1}} E^{(p,q)},\ldots,
\nabla_{\partial_{x_\ell}} E^{(p,q)}
$$
with $\deg \nabla_{\partial_{x_i}} E^{(p,q)} = p h_{1} + q h_{2}$ 
for $i = 1, \dots, \ell$, 

\end{itemize}

\end{theorem}

The existence of the {\bf primitive decomposition} of
$
D(\A,(2p-1,2q-1))^{W}
$
is proved
by the following theorem:

\begin{theorem}
\label{mainhodge}
Under the same assumption of 
Theorem \ref{basisconstruction}
define 
\begin{eqnarray*}
\theta_i^{(p,q)}:&=&
\nabla_{\partial_{P_i}}E^{(p, q)} 
= 
\nabla_{\partial_{P_i}}\nabla_{D}^{-q} \nabla_{D_1}^{q-p} E
\,\,\,\,\,\,\,\,
(1\le i\le \ell)
\end{eqnarray*}
for 
$p, q\in\Z$.
Then the set
\[
\{
\theta_{i}^{(p+k, q+k)} 
\mid
k\ge 0,
1\le i \le \ell 
\}
\]
is a $T$-basis for 
$D(\A, (2p-1, 2q-1))^{W} $.
Put
$$
\G^{(p,q)}:=\bigoplus_{i=1}^\ell T 
\cdot \theta_i^{(p,q)}.
$$
Then we have a $T$-module decomposition
(called the {primitive decomposition})
$$
D(\A,(2p-1,2q-1))^{W}
=
\bigoplus_{k \ge 0} \G^{(p+k,q+k)}.
$$
\end{theorem}

We will also prove
\begin{theorem}
For any
irreducible
 Coxeter arrangement
$\A$ and any multiplicity ${\bf m}$,
the $R$-module $D(\A,\bfm)^W$ is free.
\label{mainc}
\end{theorem}

Theorems \ref{basis} and \ref{basisconstruction} are
used to prove the freeness of Shi-Catalan arrangements
associated with any Weyl arrangements in \cite{AT10}.

The organization of this article is as follows. 
In Section \ref{proofoftheorem1.3.1}
 we prove 
Thereom \ref{basisconstruction} 
when $q \ge 0$.
In Section 3 we 
prove 
Theorem \ref{mainhodge}
to have the primitive 
decomposition, 
which is a key to complete 
the proof of Theorem 
\ref{basisconstruction}
at the end of Section 3. 
In Section 4 we verify that
the primitive derivation $D_{1}$
{can be
chosen to be}
$W$-invariant
when $\A$ is a
Coxeter arrangement of either
the type $B_\ell$ or $F_4$.
We also
review the cases of $G_{2}$ and $I_{2}(2n)\,(n\ge 4)$
and find that the primitive derivation
$D_{1} $ is $W_{2} $-antiinvariant. 
In Section 5, combining Theorem \ref{basisconstruction} with
earlier results in \cite{T02, AT08, W10},
we finally prove Theorems \ref{basis} and \ref{mainc}.

\section{Proof of Theorem \ref{basisconstruction} when $q\ge 0$ }
\label{proofoftheorem1.3.1}
In this section we prove Theorem \ref{basisconstruction} 
when $q \ge 0$. 

Recall
 $R=S^W=\R[P_1,\ldots,P_\ell]$ is the invariant ring 
with basic invariants $P_1,\ldots,P_\ell$  
such that 
$2=\deg P_1 < \deg P_2 \le \cdots \le \deg P_{\ell-1} 
< \deg P_\ell=h$, where $h$ is the Coxeter number of
$W$. Put
$D=\partial_{P_\ell} \in \Der R$
which is  a primitive derivation.
Recall
$T=\ker(D:R \rightarrow R)=
\R[P_1,\ldots,P_{\ell-1}]$. 
Then the covariant derivative $\nabla_{D} $ is $T$-linear. 
For 
$\mathbf{P}:=[P_1,\ldots,P_\ell]$, the Jacobian 
matrix $J(\mathbf{P})$ is defined as the matrix whose 
$(i,j)$-entry is $\displaystyle \frac{\partial P_j}{\partial x_i}$. 
Define 
$A:=[I^{*} (dx_i,dx_j)]_{1 \le i, j \le \ell}$ and 
$G:=[I^{*} (dP_i,dP_j)]_{1 \le i, j \le \ell}
=
J(\mathbf{P})^TAJ(\mathbf{P})$. 

\begin{define}
(\cite{Y02,W10})
Let ${\bf m} : \A\rightarrow \Z$ and  
$\zeta \in D(\A, -\infty)^{W} $.
We say that $\zeta$ is {\bf
${\bf m}$-universal}  when 
$\zeta$ is homogeneous and
the $S$-linear map
\begin{align*}
\Psi_{\zeta} : \Der&_{S} \longrightarrow D(\A, 2{\bf m})\\
&\theta \longmapsto \nabla_{\theta}\, \zeta
\end{align*}
 is bijective.
\end{define}

\begin{example}
The Euler derivation $E$
is ${\bf 0}$-universal
because ${
\Psi_{E} (\delta)} = \nabla_{\delta}\, E = \delta$
and $D(\A, {\bf 0}) = \Der_{S}$.    
\end{example}

%
%
%

Recall the $T$-automorphisms 
$$\nabla_D^k \colon D(\A,-\infty)^{W} \overset{\sim}{\longrightarrow}
 D(\A,-\infty)^{W}\ (k \in \Z)$$
from Theorem \ref{automorphism}. 
Recall the following two results concerning the $\bfm$-universality: 
\begin{theorem}
(\cite[Theorem 2.8]{W10})
\label{Wakamiko2.8}
If $\zeta$ is $\bfm$-universal, then $\nabla_{D}^{-1} \zeta $ is 
$(\bfm+{\bf 1})$-universal.  
\end{theorem}

\begin{prop}
(\cite[Proposition 2.7]{W10})
\label{Wakamiko2.7}
Suppose that $\zeta$ is $\bfm$-universal. 
Let
${\bf k} \colon 
\A \rightarrow \{+1,0,-1\}$. Then an $S$-homomorphism 
$$
\Phi_{\zeta}:
D(\A,\bfk) \rightarrow D(\A, \bfk+2\bfm)
$$
defined by 
$$
\Phi_{\zeta}(\theta):=\nabla_\theta \,
\zeta
$$
gives an $S$-module isomorphism.
\end{prop}

We require that assumption of Theorem \ref{basisconstruction}
is satisfied
in the rest of this section:
Suppose that $\A = \A_{1} \cup \A_{2} $ is the orbit decomposition
and that $D_{1}$,
a primitive derivation with respect to $\A_{1} $ 
in the sense of \cite[Definition 2.4]{AT09}, 
is $W$-invariant. 
Let
$W_i,\ R_i,\ P_{j}^{(i)},\ T_i,\ D_i\ (i=1,2)$ are 
defined as in Section \ref{introduction}. 
Even when $\A_{1} $ is not irreducible, we may consider
a $T_{1} $-isomorphism
$$\nabla_{D_{1}}^{k} \colon D(\A_{1} ,-\infty)^{W_{1}} 
\overset{\sim}{\longrightarrow}
 D(\A_{1},-\infty)^{W_{1} }\ (k \in \Z)$$
from Theorem \ref{automorphism}. 

\begin{prop}
\label{Epquniversal}
Suppose $q\ge 0$. 
The derivation $E^{(p, q)} :=
\nabla_{D}^{-q} \nabla_{D_{1}}^{q-p} E  
$ 
is $(p, q)$-universal. 
\end{prop}

\noindent
{\bf Proof.}
When $\A_{1} $ is irreducible, 
\cite{AY07} and \cite{AT08} 
 imply that 
$\nabla_{D_1}^{q-p} E 
$ is $(p-q, 0)$-universal.
When $\A_{1} $ is not irreducible, 
$\nabla_{D_1}^{q-p} E 
$ is $(p-q, 0)$-universal
because of \cite{AT09}.
Thus $E^{(p, q)} 
=
\nabla_{D}^{-q} \nabla_{D_{1}}^{q-p} E  
$ 
is $(p, q)$-universal by  
Theorem \ref{Wakamiko2.8}.
\owari

\bigskip

Since $E^{(p, q)} $ is $(p,q)$-universal, 
Proposition \ref{Wakamiko2.7}
yields the following:

\begin{prop}
Let
$q \ge 0$ and 
${\bf m} \colon 
\A \rightarrow \{+1,0,-1\}$. Then an $S$-homomorphism 
$$
\Phi_{p,q}:
D(\A,\bfm) \rightarrow D(\A, (2p, 2q)+\bfm)
$$
defined by 
$$
\Phi_{p,q}(\theta):=\nabla_\theta 
E^{(p,q)}
$$
gives an $S$-module isomorphism.
\label{main22}
\end{prop}

\noindent
{\bf Proof of Theorem \ref{basisconstruction} ($q\ge 0$).}
We may apply Proposition \ref{main22} because

(1) $\partial_{P_{1} } , \dots, \partial_{P_{\ell} } $ form a basis for 
$D(\A, (-1, -1))$,

(2) $\partial_{P^{(1)}_{1} } , \dots, \partial_{P^{(1)}_{\ell} } $
 form a basis for 
$D(\A, (-1, 0))$,

(3) $\partial_{P^{(2)}_{1} } , \dots, \partial_{P^{(2)}_{\ell} } $ form a basis for 
$D(\A, (0, -1))$, and

(4) $\partial_{x_{1} } , \dots, \partial_{x_{\ell} } $ form a basis for 
$D(\A, (0, 0))$.
\owari

\section{Primitive decompositions}

In this section we first
prove Theorem \ref{mainhodge}
to define
the primitive 
decomposition of $D(\A, (2p-1, 2q-1))^{W}$.  
Next we prove Theorem \ref{basisconstruction}.

\begin{prop}
\label{univTbasis}
Let $\zeta$ be $\bfm$-universal.  Then

(1)
the set 
$
\{\nabla_{\partial_{P_{j} } } \nabla_{D}^{-k} \zeta \mid 
1\le j\le \ell, k\ge 0
\}$
is linearly independent over
$T$.

(2)
Define
${\mathcal G}^{(k)} $ to be the free $T$-module
with basis 
$\{\nabla_{\partial_{P_{j} } } \nabla_{D}^{-k} \zeta \mid 
1\le j\le \ell\}$
for
$k\ge 0$.
Then the Poincar\'e series 
$
\Poin(\bigoplus_{k \ge 0} \G^{(k)}, t)
$ 
satisfies:
\begin{eqnarray*}
\Poin(\bigoplus_{k \ge 0} \G^{(k)}, t)
&=&
\left(
\prod_{i=1}^{\ell} \displaystyle \frac{1}{1-t^{d_i}}
\right) 
(\sum_{j=1}^\ell 
t^{p-d_{j}}),
\end{eqnarray*}
where $p=\deg \zeta$
and $d_{j} = \deg P_{j} \,\,(1\le j\le \ell)$. 

(3)
$$D(\A, 2\bfm-{1})^{W}=
\bigoplus_{k\ge 0} {\mathcal G}^{(k)}.
$$
%
\end{prop}

\noindent
\textbf{Proof.}
Let $k\in\Z_{\ge 0} $.
By Theorem \ref{Wakamiko2.8},
$\zeta^{(k)}:= \nabla_{D}^{-k} \zeta  $ 
is $(\bfm + k)$-universal, where the
``$k$'' in the
$(\bfm + k)$ stands for
the constant multiplicity $k$ by abuse of notation.
Thus by Proposition \ref{Wakamiko2.7} we have
the following two bases:
\[
\nabla_{\partial_{P_{1} } } \zeta^{(k)},
\dots
,
\nabla_{\partial_{P_{\ell} } } \zeta^{(k)}, 
\]
for the 
$S$-module 
$D(\A, 2\bfm+2k-1)$
and
\[
\nabla_{\partial_{I^{*}(dP_{1}) } } \zeta^{(k)},
\dots
,
\nabla_{\partial_{I^{*}(dP_{\ell}) } } \zeta^{(k)}, 
\]
for the $S$-module $D(\A, 2\bfm+2k+1)$.
  Note that the two bases are also $R$-bases for
$D(\A, 2\bfm+2k-1)^{W} $
and
$D(\A, 2\bfm+2k+1)^{W} $ respectively.
Since the $T$-automorphism 
$$\nabla_D \colon D(\A,-\infty)^{W} \overset{\sim}{\longrightarrow}
 D(\A,-\infty)^{W}$$
in Theorem \ref{automorphism}
induces a $T$-linear bijection
$$\nabla_D \colon D(\A, 2\bfm+2 k+1)^{W} 
\overset{\sim}{\longrightarrow}
 D(\A, 2\bfm+ 2k-1)^{W}$$
as in \cite[Theorem 4.4]{AT09},
we may find an $\ell\times \ell$-matrix
$B^{(k)}$ 
with entries in $R$ such that
\begin{align*} 
\nabla_{D} 
\left(
\left[
\nabla_{\partial_{P_{1} } } \zeta^{(k)},
\dots,
\nabla_{\partial_{P_{\ell} } } \zeta^{(k)}
\right]
G
\right)
&=
\nabla_{D} 
\left[
\nabla_{\partial_{I^{*}(dP_{1})}} \zeta^{(k)},
\dots,
\nabla_{\partial_{I^{*}(dP_{\ell})}} \zeta^{(k)}
\right]\\
&=
\left[
\nabla_{\partial_{P_{1}}} \zeta^{(k)},
\dots,
\nabla_{\partial_{P_{\ell}}} \zeta^{(k)}
\right] B^{(k)}. 
\end{align*} 
The degree of $(i,j)$-th entry of $B^{(k)}$ is $m_i+m_j-h 
\le h-2<h$. In particular, 
the degree of $B^{(k)}_{i,\ell+1-i}$ is $0$ and 
$B^{(k)}_{i,j}=0$ if $i+j < \ell+1$. 
Hence each entry of $B^{(k)}$ lies in $T$ and 
$\det B^{(k)} \in \R$.
Since $D$ is a derivation of the minimum degree
in $\Der_{R} $, one gets  $[D, \partial_{P_{i} } ] = 0$.
Thus
$\nabla_{D} \nabla_{\partial_{P_{i} } } 
=
\nabla_{\partial_{P_{i}}} \nabla_{D}$.
Operate $\nabla_{D}^{-1}$ on the both sides of the equality above,
and get
\begin{align*} 
\left[
\nabla_{\partial_{P_{1} } } \zeta^{(k)},
\dots,
\nabla_{\partial_{P_{\ell} } } \zeta^{(k)}
\right]
G
&=
\left[
\nabla_{\partial_{P_{1}}} \zeta^{(k+1)},
\dots,
\nabla_{\partial_{P_{\ell}}} \zeta^{(k+1)}
\right]
B^{(k)}. 
\end{align*} 
This implies that $\det B^{(k)}\in\R^{\times} $ 
because 
$\nabla_{\partial_{P_{1} } } \zeta^{(k)},
\dots,
\nabla_{\partial_{P_{\ell} } } \zeta^{(k)}
$  
are linearly independent over $S$.
Inductively we have  
\begin{align*} 
&~~~\left[
\nabla_{\partial_{P_{1} } } \zeta^{(k+1)},
\dots,
\nabla_{\partial_{P_{\ell} } } \zeta^{(k+1)}
\right]
=
\left[
\nabla_{\partial_{P_{1}}} \zeta^{(k)},
\dots,
\nabla_{\partial_{P_{\ell}}} \zeta^{(k)}
\right]
G
(B^{(k)})^{-1}\\
&=
\left[
\nabla_{\partial_{P_{1}}} \zeta,
\dots,
\nabla_{\partial_{P_{\ell}}} \zeta
\right]
G
(B^{(0)})^{-1}
G
(B^{(1)})^{-1}
\cdots
G
(B^{(k)})^{-1}\\
&=
\left[
\nabla_{\partial_{P_{1}}} \zeta,
\dots,
\nabla_{\partial_{P_{\ell}}} \zeta
\right]
G_{k+1},
\end{align*} 
  where 
$
G_{i} 
=
G
(B^{(0)})^{-1}
G
(B^{(1)})^{-1}
\cdots
G
(B^{(i-1)})^{-1} 
\,\,(i\ge 0)$. 
Note that $G$ appears
$i$ times in the definition of 
$G_{i} $.
For $M = (m_{ij}) \in M_{\ell}(F)$,
define $D[M] = (D(m_{ij}))
\in M_{\ell}(F)$. 
Then $D^{j} [G_{i} ] = O$ when $j > i$
and  $\det D^{i} [G_{i} ] \neq 0$ 
because $\det D[G] \neq 0$
and
$D^{2}[G]=O$ 
(e.g., see
\cite{S93, AT09}).

(1)
Suppose that 
$\{\nabla_{\partial_{P_{j} } } \zeta^{(k)} \mid 
1\le j\le \ell, k\ge 0
\}$
is linearly dependent over $T$.
Then there exist 
$\ell$-dimensional column vectors
 $
{\bfg}_{0}, 
{\bfg}_{1},
\dots
, 
{\bfg}_{q}
\in T^{\ell} 
(q\ge 0)$  
with 
$
{\bfg}_{q}
\neq {\mathbf 0}$
such that
\[
{\bf 0}
=
\sum_{i=0}^{q}  
\left[
\nabla_{\partial_{P_{1}}} \zeta^{(i)},
\dots,
\nabla_{\partial_{P_{\ell}}} \zeta^{(i)}
\right]
\bfg_{i} 
=
\left[
\nabla_{\partial_{P_{1}}} \zeta,
\dots,
\nabla_{\partial_{P_{\ell}}} \zeta
\right]
\left(
\sum_{i=0}^{q}  
G_{i} \bfg_{i}\right).
\]
Since
$\nabla_{\partial_{P_{1}}} \zeta,
\dots,
\nabla_{\partial_{P_{\ell}}} \zeta
$  
are linearly independent over $R$,
one has
$$
{\bf 0}
=
\sum_{i=0}^{q}  
G_{i} 
\bfg_{i}.
$$
 Applying the operator $D$ on the both sides
$q$ times, 
we get
$
D^{q} [G_{q}] \bfg_{q} = \bf 0.
$ 
Thus ${\mathbf g}_{q} = {\mathbf 0}$
which is a contradiction. This proves (1).
   
(2)
Compute
\begin{eqnarray*}
\mbox{Poin}(\bigoplus_{k \ge 0} \G^{(k)}, t)
&=&
\sum_{k \ge 0} 
\left(
\prod_{i=1}^{\ell-1} 
\displaystyle 
\frac{1}{1-t^{d_i}} 
\right)
(
\sum_{j=1}^\ell 
t^{p-d_{j}+ kd_{\ell}})\\
&=&
\left(
\prod_{i=1}^{\ell-1} \displaystyle \frac{1}{1-t^{d_i}}
\right)
(\sum_{k \ge 0} t^{kd_{\ell}})
 (\sum_{j=1}^\ell t^{p-d_{j}})
\\
&=&
\left(
\prod_{i=1}^{\ell} \displaystyle \frac{1}{1-t^{d_i}}
\right) 
(\sum_{j=1}^\ell 
t^{p-d_{j}}).
\end{eqnarray*}

(3)
We have
\[
D(\A, 2\bfm - 1)^{W} \supseteq \bigoplus_{k\ge 0} {\mathcal G}^{(k)} 
\]
by (1).
So it suffices to prove 
$$\mbox{Poin}(D(\A,2\bfm -1)^{W},t)
= 
\mbox{Poin}(\bigoplus_{k \ge 0} \G^{(k)},t).$$
Since $D(\A, 2\bfm-1)^{W} $ is a free $R$-module
with a basis
\[
\nabla_{\partial_{P_{1} } } \zeta,
\dots
,
\nabla_{\partial_{P_{\ell} } } \zeta, 
\]
we obtain
\begin{eqnarray*}
\mbox{Poin}(D(\A,2\bfm-1)^{W},t)&=&
\left(
\prod_{i=1}^\ell \ \displaystyle \frac{1}{1-t^{d_i}} 
\right)
(\sum_{i=1}^\ell   t^{p-d_{j}})
=
\mbox{Poin}(\bigoplus_{k \ge 0} \G^{(k)},t),
\end{eqnarray*}
which completes the proof.
\owari

\bigskip

We require that the assumption of Theorem
\ref{basisconstruction} is satisfied
in the rest of this section.

\medskip
 
\noindent
{\bf Proof of Theorem \ref{mainhodge}.}
Suppose $q\ge 0$ to begin with.
Then, 
by Proposition \ref{Epquniversal},
$E^{(p,q)}$ is 
$(p, q)$-universal.
Apply Proposition \ref{univTbasis} for
$\zeta = E^{(p,q)}$
and
$\bfm=(p, q)$, and we have Theorem 
\ref{mainhodge}:
\[
D(\A, (2p-1, 2q-1))^{W} = \bigoplus_{k\ge 0} \G^{(p+k, q+k)}
\]
when $q\ge 0$. 
Send the both handsides by $\nabla_{D} $, and we
get 
$$
D(\A, (2p-3, 2q-3))^{W} = \bigoplus_{k\ge 0} \G^{(p+k-1, q+k-1)} 
$$ 
because
$\nabla_{D} \left(D(\A, (2p-1, 2q-1))^{W}\right) = D(\A, (2p-3, 2q-3))^{W}$ 
as in \cite[Theorem 4.4]{AT09} and 
$\nabla_{D} (\theta_{i}^{(p,q)})
=
\theta_{i}^{(p-1,q-1)}  
$.
Apply $\nabla_{D}$ repeatedly to complete
the proof for all $q\in\Z$. 
\owari 

\bigskip

Note that we do not assume $p \ge 0$ in the following proposition:

\begin{prop}
For $p,q \in \Z$, 
the $S$-module  $
D(\A,(2p-1,2q-1))
$
has a $W$-invariant basis.
%
\label{dualfree}
\end{prop}

\noindent
{\bf Proof.}
Recall that 
$$
\nabla_{\partial_{P_{1}}} E^{(p,q)},
\nabla_{\partial_{P_{2}}} E^{(p,q)},
\dots
\nabla_{\partial_{P_{\ell}}} E^{(p,q)},
$$ which are $W$-invariant,
form an
$S$-basis for 
$
D(\A,(2p-1,2q-1))
$   
when $q\geq 0$ by Theorem \ref{basisconstruction} (1).
It is then
easy to see that
they are also an $R$-basis for 
$
D(\A,(2p-1,2q-1))^{W} 
$
for $q\geq 0.$
By \cite{A08} \cite{AT08}, there exists a
$W$-equivariant nondegenerate $S$-bilinear pairing
\[
(~,~)
:
D(\A,(2p-1,2q-1))
\times
D(\A,(-2p+1,-2q+1))
\longrightarrow
S,
\]
characterized   by 
\[
(I^{*}(\omega), \theta)
=
\left<
\omega, \theta
\right>
\]
where
$\omega\in
\Omega^{1} (\A, (-2p+1,-2q+1))
$ 
and 
$\theta\in
D(\A,(-2p+1,-2q+1)).
$ 
Let 
$\theta_{1}, \dots , \theta_{\ell}  $
 denote 
the dual basis for 
$D(\A,(-2p+1,-2q+1))$ 
satisfying
\[
\left(
\nabla_{\partial_{P_{i}}} E^{(p,q)},
\theta_{j} 
\right)
=
\delta_{ij} 
\]
for $1\leq i, j \leq \ell$.
Then
$\theta_{1}, \dots , \theta_{\ell}  $
are $W$-invariant because the pairing
$(~,~)$ is $W$-equivariant.   
%
\owari 

\bigskip

Although the following lemma is standard and easy,
we give a proof for completeness.
 
\begin{lemma}
\label{MWfree} 
Let $M$ be an S-submodule of
$\Der_{F}$.  The following two conditions are equivalent:

(1) $M$ has a $W$-invariant basis $\Theta$ over $S$.   

(2) The $W$-invariant part $M^{W} $ is a free $R$-module with a basis $\Theta$
and there exists a natural $S$-linear isomorphism
\[
M^{W} \otimes_{R} S \simeq 
M.
\]
\end{lemma}

\medskip

\noindent
{\bf Proof.}
It suffices to prove that (1) implies (2)
because the other implication is obvious.
Suppose that 
$\Theta =\{\theta_{\lambda}\}_{\lambda\in\Lambda}
$ is a $W$-invariant basis for 
$M$ over $S$.
Since it is linearly independent over $S$,
so is over $R$.     
Let $\theta\in M^{W}$.
Express
\[
\theta = \sum_{i=1}^{n} f_{i} \theta_{i}  
\]
 with $f_{i} \in S$ and $ \theta_{i} \in\Theta\,\,(i=1,\dots,n)$.
 Let $w\in W$ act on the both handsides.
Then we get
\[
\theta = \sum_{i=1}^{n} w(f_{i}) \theta_{i}.  
\]
This implies $f_{i} = w(f_{i} )$ for every $w\in W$.
Hence
$f_{i} \in R$ for each $i$.  
Therefore
$\Theta$
is a basis for $M^{W} $ over $R.$
This is (2).
\owari

\bigskip

\begin{prop}
\label{Epquniversal} 
For any  $p, q\in\Z$, $E^{(p, q)} $ is $(p, q)$-universal.  
\end{prop}

\noindent
{\bf Proof.}
By Theorem \ref{mainhodge} we have the decomposition:
\[
D(\A, (2p-1, 2q-1))^{W} 
= 
\bigoplus_{k\ge 0} \G^{(p+k, q+k)} 
\]
for $p, q\in\Z$. 
As we saw in Proposition \ref{univTbasis}
(2),
we have
\begin{multline}
\label{3.1} 
\mbox{Poin}(D(\A, (2p-1, 2q-1))^{W}, t)
=
\mbox{Poin}(\bigoplus_{k \ge 0} \G^{(p+k, q+k)}, t)\\
=
\left(
\prod_{i=1}^{\ell} \displaystyle \frac{1}{1-t^{d_i}}
\right) 
(\sum_{i=1}^\ell 
t^{m- d_{j}}),
\end{multline} 
where $m := \deg E^{(p,q)}.$ 
Recall that the $S$-module 
$
D(\A, (2p-1, 2q-1))
$
has a $W$-invariant basis
$\theta_{1}, \dots, \theta_{\ell}$
by
Proposition \ref{dualfree}.
By Lemma \ref{MWfree}, we know that
$\theta_{1}, \dots, \theta_{\ell}$
form a basis 
 for the $R$-module 
$
D(\A, (2p-1, 2q-1))^{W}
$.
Thanks to 
(\ref{3.1})
%
we may assume
that
$
\deg \theta_{j} = 
m-d_{j} 
=
\deg \nabla_{\partial_{P_{j} } } E^{(p,q)}. 
$ 
Therefore there exists $M \in M_\ell(R)$ such that 
$$
[\theta_1,\ldots,\theta_\ell]M
=[
\nabla_{\partial_{P_{1} } } E^{(p,q)},
\ldots,
\nabla_{\partial_{P_{\ell} } } E^{(p,q)}
]
$$
with $\det M\in\R$. 
Since
$$
\max_{1 \le i, j \le \ell} \left|
\deg \theta_{i}
-
\deg \nabla_{\partial_{P_{j} } } E^{(p, q)}\right|
=d_{\ell} - d_{1} < \deg P_{\ell},
$$
we get $M\in M_{\ell} (T)$. 
Since $
\nabla_{\partial_{P_{1} } } E^{(p,q)},
\ldots,
\nabla_{\partial_{P_{\ell} } } E^{(p,q)}
$ 
are
linearly independent over $T$
by Proposition \ref{univTbasis} (1), 
we have
 $\det M\in \R^{\times}$. 
Thus
$$\nabla_{\partial_{P_{1} } } E^{(p,q)},
\ldots,
\nabla_{\partial_{P_{\ell} } } E^{(p,q)}
$$
form an $S$-basis for $
D(\A, (2p-1, 2q-1)).$
Since
\[
\left[
\nabla_{\partial_{P_{1} } } E^{(p,q)},
\ldots,
\nabla_{\partial_{P_{\ell} } } E^{(p,q)}
\right]
J({\mathbf P})^{T} 
=
\left[
\nabla_{\partial_{x_{1} } } E^{(p,q)},
\ldots,
\nabla_{\partial_{x_{\ell} } } E^{(p,q)}
\right],
\]
we may
apply
the multi-arrangement version of 
  Saito's criterion
\cite{S80, Z89, A08}
to prove that
$\nabla_{\partial_{x_{1} } } E^{(p,q)},
\ldots,
\nabla_{\partial_{x_{\ell} } } E^{(p,q)}
$
form an $S$-basis for $
D(\A, (2p, 2q))$
for any $p, q\in\Z$.
This shows that 
$E^{(p,q)} $ is $(p, q)$-universal
for any $p, q\in\Z$.
 \owari

\bigskip
\noindent
\textbf{Proof of Theorem \ref{basisconstruction}}
($q\in\Z$).
Theorem \ref{Wakamiko2.8} and  Proposition \ref{Epquniversal} 
complete the proof by the same argument as 
that in Section 2 for $q\ge 0$. \owari
\medskip

\section{The cases of $B_\ell$,  $F_{4}$, $G_{2} $ 
and 
$I_{2} (2n)$ }
\label{invariantinvariant1}
\noindent
$\bullet$ {\bf The case of $B_{\ell} $}

\bigskip

The roots of the type $B_{\ell}$
are:
\[
\pm x_{i}, \,
\pm x_{i} \pm x_{j}\,
\,\,\,\,
(1\le i < j \le \ell)
\]
in terms of an orthonormal basis
$x_{1} ,  \dots, x_{\ell}  $  for $V^{*} $. 
Altogether there are $2 \ell^{2} $  of them.
Define
\[
Q_1:=
\prod_{i=1}^{\ell}  x_i ,
\,\,\,
Q_2:=\prod_{1\le i<j\le\ell}(x_i \pm x_{j}),
\,\,\,
Q = Q_{1} Q_{2}. 
\]
Then the arrangement $\A_{1} $ defined by $Q_{1} $ is
of the type $A_{1} \times\dots\times A_{1} = A_{1}^{\ell}$.
The arrangement $\A_{2} $ defined by $Q_{2} $ is
of the type $D_{\ell}$.
The arrangement
$\A$
defined by $Q$ is 
of the type $B_{\ell} $ 
and
$\A = \A_{1} \cup \A_{2} $ 
is the orbit decomposition.
Note that $A_{1}^{\ell}  $ 
is not irreducible.
{Define
$$
D_1:=\sum_{i=1}^\ell \displaystyle \frac{1}{x_i}\partial_{x_i}
$$
which} is a primitive derivation in the sense of 
\cite{AT09}. 
Obviously
$D_1$ is $W$-invariant. 
Let
$P_{j} = \sum_{i=1}^{\ell} x_{i}^{2j}  
\,\,(j\ge 1) $. 
 Then
$P_1,\ldots,P_{\ell}
$ form a set of basic invariants under $W$
while
$Q_1,P_1,\ldots,P_{\ell-1}
$ form a set of basic invariants under $W_{2} $.
Define a primitive derivation $D_{2} $ with respect to
$\A_{2} $ so that
$$
D_{2} (Q_{1} )
=
D_{2} (P_{j})
=
0
\,\,(j=1,\dots,\ell-2),
\,\,\,
D_{2} (P_{\ell-1}
)=1.$$ 
Thus
$$(w D_{2}) (P_{\ell-1}
)
=
D_{2} (w^{-1}  P_{\ell-1}
)
=
D_{2} (P_{\ell-1})
=1
\,\,\,(w\in W).$$ 
This implies that $D_{2} $ is $W$-invariant.

\bigskip
\noindent
$\bullet$ {\bf The case of $F_{4} $}

\bigskip

The roots of the type $F_{4}$
are:
\[
\pm x_{i}, \,
(\pm x_{1}\pm x_{2}\pm x_{3}\pm x_{4})/2,\,
\pm x_{i} \pm x_{j}\,
\,\,\,\,(1\le i < j \le 4)
\]
in terms of an orthonormal basis
$x_{1} , x_{2} , x_{3}, x_{4}  $  for $V^{*} $. 
Altogether there are 48 of them.
Define
\[
Q_1:=
\prod_{1 \le i < j  \le 4} (x_i \pm x_j),
\,
Q_2:=\prod_{i=1}^4 x_i \prod(x_1 \pm x_2 \pm x_3 \pm x_4),
\,
Q = Q_{1} Q_{2}. 
\]
The arrangement $\A_{i} $ defined by $Q_{i} $ is
of the type $D_{4}$ 
$(i=1,2)$.
Then the arrangement
$\A$
defined by $Q$ is 
of the type $F_{4} $ 
and
$\A = \A_{1} \cup \A_{2} $ 
is the orbit decomposition.
Define
\begin{eqnarray*}
P_1^{(1)}&=&
\sum_{i=1}^{4}  x_i^2,\,\,
P_2^{(1)}=
\sum_{i=1}^{4}  x_i^4,\,\,
P_3^{(1)}=x_{1} x_{2} x_{3} x_{4},\,\,
P_4^{(1)}=
\sum_{i=1}^{4} x_{i}^{6} 
+ 5
\sum_{i\neq j} x_{i}^{2} x_{j}^{4}.     
\end{eqnarray*}
Compute
\[
P_{4}^{(1)} 
=
-4
\sum_{i=1}^{4} x_{i}^{6} 
+5
P_1^{(1)}P_2^{(1)}.
\]
Thus 
$P_1^{(1)},
P_2^{(1)},
P_3^{(1)},
P_{4}^{(1)}
$ 
are a set of basic invariants under $W_{1}$.
The reflection $\tau$ with respect to $x_1+x_2+x_3+x_4=0$ is given by 
$$
\tau(x_i)=\displaystyle \frac{2x_i-\sum_{j=1}^4 x_j}{2}\ (i=1,2,3,4).
$$
A calculation shows that $P_{4}^{(1)}$ 
is $\tau$-invariant. 
Let $s_{i}$ denote 
the reflection 
with respect to $x_{i} = 0\,\,(1\le i\le 4)$.
Since the Coxeter group $W_{2} $ is generated by
$\tau $ and $s_{i} \,\,(1\le i\le 4)$,
we know that $P_{4}^{(1)}$ is $W_{2}$-invariant thus
$W$-invariant.   
Define a primitive derivation $D_{1} $ with respect to 
$\A_{1}$ so that
$$D_{1} (P_{j}^{(1)}  )
=
0
\,\,(j=1,2,3),
\,\,\,
D_{1} (P_{4}^{(1)}
)=1.$$ 
Thus
$$(w D_{1}) (P_{4}^{(1)}
)
=
D_{1} (w^{-1}  P_{4}^{(1)}
)
=
D_{1} (P_{4}^{(1)}
)
=1
\,\,\,(w\in W).$$ 
This implies that $D_{1} $ is $W$-invariant.
We conclude that $D_{2}$ is also $W$-invariant because
an orthonormal coordinate change
\[
x_{1} = \frac{y_{1}-y_{2}}{\sqrt{2}},
\,\,
x_{2} = 
\frac{y_{1}+y_{2}}{\sqrt{2}},
\,\,
x_{3} = \frac{y_{3}-y_{4}}{\sqrt{2}},
\,\,
x_{4} = \frac{y_{3}+y_{4}}{\sqrt{2}}
\]
 switches $\A_{1} $ and $\A_{2}.$

\bigskip
\noindent
$\bullet$ {\bf The cases of $G_{2} $ and $I_{2}(2n)\, (n\ge 4)$}

\bigskip

The arrangement $\A$ of the type $G_{2} $
has exactly two orbits
$\A_{1} $ 
and
$\A_{2} $, 
each of which is of the type $A_{2} $.
Let $n\ge 4$.
Then the arrangement $\A$ of the type $I_{2} (2n)$ 
has exactly two orbits
$\A_{1} $ 
and
$\A_{2} $, 
each of which is of the type 
$I_{2} (n)$. 
In both cases,
by \cite{W10}, one may choose
\[
D_{1} = Q_{2} D,
\,\,\,
D_{2} = Q_{1} D.
\]
 Since $Q_{2} $ is $W_{2} $-antiinvariant and $D$ is $W$-invariant,
$D_{1} $ is $W_{2} $-antiinvariant.   
Similarly
$D_{2} $ is $W_{1} $-antiinvariant.

\section{Proofs of Theorems \ref{basis} and \ref{mainc} }

Assume that $\A$ is an irreducible
Coxeter arrangement
in the rest of the article.

\medskip

\noindent
\textbf{Proof of Theorem \ref{basis}}. 
If $\A$ has the single orbit,
%
then 
the result in
 \cite{T02, AY07, AT08} 
completes the proof. If not, then 
$\A$ has exactly two orbits.
If $\A$ is of the type either $G_{2} $ or
$I_{2}(2n)$ with $n\ge 4$, 
then 
$D(\A, \bfm)$ is a free $S$-module because
$\A$ lies in a two-dimensional vector space. 
For the remaining cases of the type  $B_{\ell} $ and $F_{4} $, 
Section 4 allows us to apply
Theorem \ref{basisconstruction} 
to complete the proof. \owari

\bigskip

A multiplicity $\bfm : \A \rightarrow \Z$ is said to be 
{\bf odd} 
if its image lies in $1+2\Z$.
 
\begin{prop}
\label{equivariantodd} 
If $\bfm$ is  equivariant and odd, then 
$D(\A, \bfm)$ has a $W$-invariant basis
over $S$.
\end{prop}

\medskip

\noindent
{\bf Proof.}
When $\A$ has the single orbit, $\bfm$ is constant.
In this case
Proposition was proved in
\cite{T02, AY07, AT08}.
If $\A$ is of the type
either $G_{2} $ or $I_{2}(2n)\,(n\ge 4)$,
then
Proposition was verified in \cite{W10}. 
For the remaining cases of 
 $B_{\ell} $ and
 $F_{4} $, 
Proposition \ref{dualfree} 
completes the proof. 
\owari

\bigskip

Recall the $W$-action on $\A$:
\[
W\times\A \longrightarrow \A
\]
by sending $(w, H)$ to $wH$ ($w\in W$, $H\in\A$).
For any multiplicity $\bfm: \A \rightarrow \Z$, define 
a new multiplicity $\bfm^{*}$ by
\[
\bfm^{*} (H)
:=
\max_{w\in W} \left(
2\cdot \lfloor \bfm(wH)/2 \rfloor
+1
\right),
\]
where $\lfloor a \rfloor$ stands for the greatest integer 
not exceeding $a$. 
Then $\bfm^{*} $  is obviously equivariant and odd.

\begin{prop}
\label{mstar}
For any irreducible Coxeter arrangement $\A$ and any multiplicity $\bfm$,
$$
D(\A, \bfm)^{W} 
=
D(\A, \bfm^{*} )^{W}. 
$$
\end{prop}

\noindent
{\bf Proof.}
Since $\bfm(H) \le \bfm^{*} (H)$ for any $H\in\A$,
we have $$D(\A, \bfm)^{W}  \supseteq D(\A, \bfm^{*})^{W}.$$
We will show the other inclusion.
Let $H\in\A$ and $\theta\in D(\A, \bfm)^{W}$.
It suffices to verify the following two statements:

\medskip

(A) $\theta(\alpha_{H}) \in \alpha_{H}^{\bfm(wH)} S_{(\alpha_{H} )}  $ 
for any $w\in W$,

(B) $\theta(\alpha_{H}) \in \alpha_{H}^{2m} S_{(\alpha_{H} )}  $ implies
$\theta(\alpha_{H}) \in \alpha_{H}^{2m+1} S_{(\alpha_{H} )}  $ 
for any $m\in\Z$. 
 
\medskip

For $w\in W$ let $w^{-1}$ act on the both sides of
$$\theta(\alpha_{wH}) \in \alpha_{wH}^{\bfm(wH)} S_{(\alpha_{wH} )}  $$ 
to get
$$\theta(\alpha_{H}) \in \alpha_{H}^{\bfm(wH)} S_{(\alpha_{H} )}.  $$ 
This verifies (A).

Fix $H\in\A.$ 
Let $s$ be the orthogonal
reflection through $H$. Then $s(\alpha_{H}) = -\alpha_{H}$ . 
Suppose that 
$\theta(\alpha_{H}) = \alpha_{H}^{2m} p  $
with $p\in S_{(\alpha_{H})}$. 
Let $s$ 
act on the both handsides and we have
$\theta(-\alpha_{H}) = (-\alpha_{H})^{2m} s(p)  $.
This implies $-p=s(p)$.   Since $s(p) = p$ on $H$,
one has $p=0$ on $H$,
which implies
$p\in \alpha_{H} S_{(\alpha_{H})}$. 
This verifies (B).
\owari

\medskip

\noindent
{\bf
Proof of Theorem \ref{mainc}. 
}
Thanks to Proposition \ref{mstar} 
we may assume
that
$\bfm$ is equivariant and
odd.
Apply Proposition \ref{equivariantodd}
and Lemma \ref{MWfree}.
\owari

\medskip

\begin{cor}
\[
D(\A, \bfm)^{W} \otimes_{R} S
\simeq 
D(\A, \bfm^{*}).
\]
\label{cor2} 
\end{cor}

\medskip

\noindent
{\bf
Proof. 
}
Apply Proposition \ref{equivariantodd}
and Lemma \ref{MWfree} to get
\[
D(\A, \bfm^{*})^{W} \otimes_{R} S
\simeq 
D(\A, \bfm^{*}).
\]
Then Proposition \ref{mstar} 
completes the proof.
\owari

\bigskip

The following corollary shows that
the converse of Proposition \ref{equivariantodd} is true.

\begin{cor}
The $S$-module $
D(\A, \bfm)$ has a $W$-invariant basis if and only if
$\bfm$ is odd and equivariant.
\end{cor}

\medskip

\noindent
{\bf
Proof. 
}
Assume that
$
D(\A, \bfm)$ has a $W$-invariant basis
over $S$.
Then, 
by Lemma \ref{MWfree}, 
we get
\[
D(\A, \bfm)^{W} \otimes_{R} S
\simeq 
D(\A, \bfm).
\]
Compare this
with Corollary \ref{cor2} and we 
know that
there exsits a common $S$-basis for 
both
$D(\A, \bfm)$ and $ D(\A, \bfm^{*})$.
By the multi-arrangement version of Saito's criterion
\cite{S80, Z89, A08},
we have $\bfm=\bfm^{*}$.
\owari

 \vspace{5mm}

\end{document}